\documentclass[11pt, a4paper]{article}

\usepackage{amssymb, amsmath, textcomp, amsthm}
\usepackage{bbm,dsfont}
\usepackage[pdftex]{graphicx}
\numberwithin{equation}{section}

\title{Remarks on the maximal operator and Hilbert transform along variable parabolas}
\author{Shaoming Guo}
\date{}

\def\R{\mathbb{R}}
\def\N{\mathbb{N}}

\def\Z{\mathbb{Z}}

\def\lesim{\lesssim}

\def\beq{\begin{equation}}
\def\endeq{\end{equation}}
\def\mc{\mathcal}

\theoremstyle{plain}
\newtheorem{thm}{Theorem}[section]

\newtheorem{lem}[thm]{Lemma}

\newtheorem{rem}[thm]{Remark}

\newtheorem{question}[thm]{Question}
\newtheorem*{openproblem*}{Open Problem}

\begin{document}
\maketitle

\begin{abstract}
We prove the boundedness of the maximal operator and Hilbert transform along certain variable parabolas in $L^p$ for $p>p_0$ with some $p_0\in (1, 2)$. Connections with the Hilbert transform along vector fields and the polynomial Carleson's maximal operator are also discussed.
\end{abstract}

\let\thefootnote\relax\footnote{Date: \date{\today}}

\section{Introduction and statement of the main result}
On the plane $\R^2$, let $u:\R^2\to \R$ be a measurable function, let $\gamma:\R\to \R$ with $\gamma(0)=0$ be a convex function, the main objects in the present paper are the following maximal operator and Hilbert transform along the variable curve $(t, u\cdot \gamma(t))$, defined separately by
\beq\label{0504ee1.1}
\mc{M}^u_{\gamma}f(x, y):=\sup_{k\in \Z} \frac{1}{2^k} \left| \int_{-2^k}^{2^k} f(x-t, y-u(x, y)\cdot \gamma(t))dt\right|,
\endeq
and
\beq\label{0504ee1.2}
\mc{H}^u_{\gamma}f(x, y):=\int_{\R}f(x-t, y-u(x, y)\cdot \gamma(t))\frac{dt}{t}.
\endeq
The case $u\equiv 1$ has been intensively studied by many authors, see for example \cite{CCC}, \cite{NRW1}, \cite{NRW2}, \cite{NVWW1} and \cite{NVWW2}. This is a translation-invariant case. Concerning the non-translation-invariant case: when $u(x, y)=x$, the operators \eqref{0504ee1.1} and \eqref{0504ee1.2} were considered by Carbery, Wainger and Wright \cite{CWW2}, which was motivated by their previous work \cite{CWW1} on the Heisenberg group; the case $u(x, y)=P(x)$ where $P:\R\to \R$ is a polynomial was studied by Bennett \cite{Bennett}. For more related work, see \cite{CP}, \cite{Kim}, \cite{Seeger} and \cite{SeegerWainger}.\\

Our intention is to study the case when $u$ is a rough function instead of a polynomial (compared with the results in \cite{Bennett} and \cite{CWW2}). To do this, we choose the simplest setting where
\beq
\gamma(t)=t^2,
\endeq
and consider 
\beq\label{0504ee1.4}
\mc{M}^uf(x, y):=\sup_{k\in \Z} \frac{1}{2^k} \left| \int_{-2^k}^{2^k} f(x-t, y-u(x, y)\cdot t^2)dt\right|,
\endeq
and
\beq\label{0504ee1.5}
\mc{H}^uf(x, y):=\int_{\R}f(x-t, y-u(x, y)\cdot t^2)\frac{dt}{t}.
\endeq
For certain measurable functions $u:\R^2\to \R$, by taking the function $f$ to be the indicator function of the unit ball $B_1(0)$, it is not difficult to see that $\mc{M}^u f$ and $\mc{H}^u f$ might not lie in $L^p$ for any $p>1$. Hence except for just measurability, we need to add more restrictions on the function $u$. 

The restriction we impose here is that $u:\R^2\to \R$ is constant in its second variable $y$. To state our results, we first introduce some notations. Let $\phi_0:\R\to \R$ be a smooth bump function supported on $[-5/2, -1/2]\cup [1/2, 5/2]$ with
\beq
\phi_0(t)=1, \forall t\in [-2, -1]\cup [1, 2].
\endeq 
For $k\in \Z$, denote
\beq\label{0504ee1.7}
\phi_k(t)=\frac{1}{2^k} \phi_0(\frac{t}{2^k}).
\endeq
Moreover, let $P_k$ denote the Littlewood-Paley projection operator in the second variable, i.e. 
\beq
P_k f(x, y):= \int_{\R} f(x, y-t)\check{\phi}_k (t)dt.
\endeq
\begin{thm}\label{2604theorem1.1}
Let $u: x\in \R\to u(x)\in \R$ be an arbitrary measurable function. Let $\mc{M}^u$ and $\mc{H}^u$ be given by \eqref{0504ee1.4} and \eqref{0504ee1.5} separately. Then there exists $p_0\in (1, 2)$ such that for all $p>p_0$, there exists a constant $C_{p}>0$ depending only on $p$ such that 
\beq\label{0402ee1.1}
\left\| \mc{M}^u \circ P_k (f) \right\|_{p} \le C_p \|f\|_p,
\endeq
and 
\beq\label{0504ee1.12}
\left\| \mc{H}^u \circ P_k (f) \right\|_p \le C_p \|f\|_p,
\endeq
for all $k\in \Z$.
\end{thm}

One motivation of studying the above special case of the function $u$ is from the following results due to Bateman \cite{Bateman} and Bateman and Thiele \cite{BT} on the boundedness of the Hilbert transform along the one-variable vector fields. For a given measurable function $u:\R\to \R$, define the Hilbert transform along the one-variable vector field $(1, u)$ by
\beq
H^u f(x, y):=\int_{\R}f(x-t, y-u(x)t)\frac{dt}{t}.
\endeq
Building on Lacey and Li's work \cite{LL0} and \cite{LL}, Bateman \cite{Bateman} proved
\begin{thm}(\cite{Bateman})\label{2604theorem1.2}
For any $p>1$, for any measurable function $u:\R\to \R$, there exists a positive constant $C_p$ depending only on $p$ such that
\beq
\|H^u \circ P_k (f)\|_p \le C_p \|f\|_p,
\endeq
for all $k\in \Z$.
\end{thm}
Basing on Theorem \ref{2604theorem1.2} and a new square function estimate, Bateman and Thiele \cite{BT} proved 
\begin{thm}(\cite{BT})\label{batemanthiele}
For any $p>3/2$, for any measurable function $u:\R\to \R$, there exists a positive constant $C_p$ depending only on $p$ such that
\beq
\|H^u f\|_p \le C_p \|f\|_p.
\endeq
\end{thm}
Hence our Theorem \ref{2604theorem1.1} is a partial generalisation of Bateman's result in \cite{Bateman} from the one-variable vector field $(t, u(x)t)$ to the ``one-variable parabolas'' $(t, u(x)t^2)$.
\begin{rem}
Unlike Bateman's result \cite{Bateman}, in Theorem \ref{2604theorem1.1} we do know how to prove the $L^p$ boundedness for all $p>1$. However our results do include the strong type $L^2$ boundedness. Indeed the $L^2$ boundedness is the main obstacle in Bateman's generalisation of Lacey and Li's results \cite{LL0} and \cite{LL}, which have already contained the weak type $L^2$ boundedness. And It is exactly with this strong $L^2$ boundedness that Bateman and Thiele \cite{BT} managed to organise all the frequency annuli together to obtain Theorem \ref{batemanthiele}. So far we do not know how to adapt Bateman and Thiele's argument to get rid of the frequency restriction in the estimates \eqref{0402ee1.1} and \eqref{0504ee1.12}.
\end{rem}

Another motivation of studying the operators \eqref{0504ee1.4} and \eqref{0504ee1.5} with $u$ being constant in $y$ is from their connection with Stein and Wainger's polynomial Carleson's maximal operator in \cite{SteinWainger}. One case that Stein and Wainger \cite{SteinWainger} proved is that for any $p>1$, there exists $C_p>0$ depending only on $p$ such that 
\beq\label{0504ee1.9}
\left\|\sup_{A\in \R}\left|\int_{\R}f(x-t)e^{iAt^2}\frac{dt}{t} \right|\right\|_{L^p(\R)} \le  C_p \|f\|_{L^p(\R)}.
\endeq
By Plancherel's theorem, the $L^2$ estimate in \eqref{0504ee1.9} is equivalent to the $L^2$ estimate of the operator $\mc{H}^u\circ P_k$ in \eqref{0504ee1.12}. We refer to \cite{Guo} for the detailed discussion. Due to the use of Plancherel's theorem, this equivalence holds only at the $L^2$ level. In Theorem \ref{2604theorem1.1}, we extend the $L^2$ bound of $\mc{H}^u\circ P_k$ to $L^p$ for all $p$ greater than certain $p_0\in (1, 2)$.\\

In the definition of the maximal operator $\mc{M}^u$, we are taking averages along parabolas. Locally this is essentially the same as taking averages along circles of proper radii. Hence the estimates in Theorem \ref{2604theorem1.1} have close connections with the boundedness of the circular maximal function (see Stein \cite{Stein2}, Bourgain \cite{Bourgain2} and Schlag \cite{Schlag}) and the Radon variant of the Kakeya problem (see Kolasa and Wolff \cite{KW} and Wolff \cite{Wolff2}). This serves as a third motivation of studying the operators \eqref{0504ee1.4} and \eqref{0504ee1.5}. \\

{\bf Notations:} Throughout this paper, we will write $x\lesim y$ to mean that there exists a constant $C$ such that $x\le C y$, and $x\sim y$ to mean that $x\lesim y$ and $y\lesim x$.  $\mathbbm{1}_E$ will always denote the characteristic function of the set $E$.\\

\section{Maximal operator case: Proof of \eqref{0402ee1.1}}

In this section, we will prove the estimate \eqref{0402ee1.1} in Theorem \ref{2604theorem1.1}. That is, for an arbitrary measurable function $u:\R\to \R$, define 
\beq\label{0604ee2.1}
\mc{M}^u f(x, y):=\sup_{k\in \Z} \frac{1}{2^k} \left| \int_{-2^k}^{2^k} f(x-t, y-u(x)\cdot t^2)dt\right|,
\endeq
we will show that there exists a $p_0<2$ such that for all $p>p_0$, it holds that
\beq\label{0504ee2.2}
\left\|\mc{M}^u \circ P_k (f)\right\|_p \lesim \|f\|_p,
\endeq
where the constant is independent of $k\in \Z$. \\
%

By the anisotropic scaling
\beq\label{0504ee2.3}
x\to x, y\to \lambda y,
\endeq
it suffices to prove \eqref{0504ee2.2} for $k=0$. In the rest of this section, we make the convention that whenever we use $f$ to denote a function, it always holds that 
\beq
P_0 f=f.
\endeq
Moreover, recall the definition of the $L^1$ normalised bump function $\phi_k$ in \eqref{0504ee1.7}, by losing a constant factor, it is not difficult to see that $\mc{M}^u$ can be bounded by
\beq\label{1602ee1.1}
\sup_{k\in \Z}\left|\int_{\R} f(x-t, y-u(x)t^2)\phi_k(t)dt\right|.
\endeq
For the sake of simplicity, we will still use $\mc{M}^u f$ to denote \eqref{1602ee1.1}. Hence, the estimate \eqref{0402ee1.1} is reduced to 
\beq
\|\mc{M}^u f\|_p \lesim \|f\|_p.
\endeq

%
%
%
%
%
%
%
%
%
%
%
\noindent For fixed $u\in \R$ and $k\in \Z$, define
\beq\label{1602ee1.3}
\mc{M}^u_k f(x, y):=\left|\int_{\R}f(x-t, y-u\cdot t^2)\phi_k(t)dt\right|.
\endeq
After a linearisation of the maximal operator $\mc{M}^u$ in \eqref{1602ee1.1}, we are to prove a bound for
\beq\label{1602ee1.4}
\mc{M}^{u(x)}_{k(x, y)}f(x, y),
\endeq
where $k:\R^2\to \Z$ denotes an arbitrary measurable function.\\

Before getting into the details of the proof, let us first explain the main ideas. Depending on the values of the functions $u$ and $k$, there will be two cases: 
\beq
u(x)\cdot 2^{2k(x, y)}\le 1 \text{ or } u(x)\cdot 2^{2k(x, y)}\ge 1,
\endeq
which will be called Case 1 and Case 2 separately.

The idea to handle Case 1 is that for 
\beq
t\in [2^{k(x, y)}, 2^{k(x, y)+1}],
\endeq
the function $u(x)t^2$ will always be ``small''. Hence we will compare the integration 
\beq
\int_{\R}f(x-t, y-u(x) t^2)\phi_{k(x, y)}(t)dt
\endeq
with 
\beq\label{1602ee1.8}
\int_{\R}f(x-t, y)\phi_{k(x, y)}(t)dt,
\endeq 
by taking the advantage that the function $f$ has frequency one in the second variable. Moreover, the expression \eqref{1602ee1.8} can be bounded by the strong maximal operator $M_S$ on the plane.

The idea to handle Case 2 is to use techniques from the oscillatory integrals: if we first freeze the linearisation functions $u$ and $k$, the expression \eqref{1602ee1.3} is of a convolution form, and its multiplier is given by
\beq
m^u_k(\xi, \eta):=\int_{\R} e^{it\xi+i u\cdot t^2 \eta}\phi_k(t)dt=\int_{\R} e^{2^k t\xi+ i u\cdot 2^{2k} t^2\eta}\phi_0(t)dt.
\endeq
By applying Van der Corput's lemma (see Page 332 in Stein \cite{Stein}), the multiplier $m^u_k$ will have certain exponential decay.\\

We proceed with the estimate of \eqref{1602ee1.4}. By losing a factor of two, it suffices to prove the boundedness of \eqref{1602ee1.4} under the assumptions 
\beq\label{0402ee1.5}
u(x)\cdot 2^{2k(x, y)}\le 1, \text{ a. e. } \R^2,
\endeq
and 
\beq\label{0402ee1.6}
u(x)\cdot 2^{2k(x, y)}\ge 1, \text{ a. e. } \R^2
\endeq
separately.

\subsection{Case \eqref{0402ee1.5}}

Following the idea explained above, we compare $\mc{M}^{u(x)}_{k(x, y)}$ with 
\beq\label{0402ee1.7}
\int_{\R}f(x-t, y)\phi_{k(x, y)}(t)dt,
\endeq
as we are in the case that the ``perturbation'' $u(x)\cdot 2^{2k(x, y)}$ to the above expression is ``small''.  
\beq\label{0402ee1.8}
\begin{split}
& \left|\int_{\R}f(x-t, y-u(x)t^2)\phi_{k(x, y)}(t)dt-\int_{\R}f(x-t, y)\phi_{k(x, y)}(t)dt\right|\\
& \le \left|\int_{\R}\int_{\R}f(x-t, \eta)\left( \check{\phi}_0(y-u(x)t^2-\eta)-\check{\phi}_0(y-\eta) \right) \phi_{k(x, y)}(t)dtd\eta\right|,
\end{split}
\endeq
where we have used the fact that $P_0 f=f$. By the non-stationary phase method, it is not difficult to see that 
\beq\label{0504ee2.18}
\left| \check{\phi}_0(y-u(x)t^2-\eta)-\check{\phi}_0(y-\eta) \right|\lesim \sum_{h\in \N} \frac{1}{(h+1)^N} \mathbbm{1}_{[y-h-1, y+h+1]}(\eta),
\endeq
where $N\in \N$ is some large integer. 
Hence the right hand side of \eqref{0402ee1.8} can be bounded by the strong maximal operator $M_S$, which satisfies the $L^p$ boundedness for all $p>1$. Moreover, the term \eqref{0402ee1.7} can also be bounded by the strong maximal operator. So far we have finished the proof of Case \eqref{0402ee1.5}.

\subsection{Case \eqref{0402ee1.6}}
In this case, we will apply the Fourier method. To do this, we first replace the $l^{\infty}$ norm over $k\in \Z$ in the definition of the maximal operator $\mc{M}^u$ by an $l^2$ norm, namely
\beq
\begin{split}
\mc{M}^{u(x)}_{k(x, y)}& = \left|\frac{1}{2^{k(x, y)}}\int_{\R}f(x-t, y-u(x)t^2)\psi_{k(x, y)}(t)dt\right|\\
	& \lesim \left( \sum_{l\in \N}\left|\int_{\R}f(x-t, y-u(x)t^2)\psi_l(u(x)^{1/2}t)\frac{dt}{|t|}\right|^2 \right)^{1/2},
\end{split}
\endeq
where $\psi_l(t):=\phi_0(t/2^l)$.  By a rather brutal estimate, we observe that 
\beq\label{0504ee1.19}
\left|\int_{\R}f(x-t, y-u(x)t^2)\psi_l(u(x)^{1/2}t)\frac{dt}{|t|}\right| \lesim 2^{2l}\cdot M_S f(x, y).
\endeq
Hence to prove the $L^p$ bounds of $\mc{M}^{u(x)}_{k(x, y)}$ for $p>p_0$ with certain $p_0<2$, by interpolation and the triangle inequality, 
it suffices to prove that there exists a universal constant $\gamma>0$ such that
\beq\label{1602ff1.20}
\left\|\int_{\R}f(x-t, y-u(x)t^2)\psi_l(u(x)^{1/2}t)\frac{dt}{|t|}\right\|_2 \lesim 2^{-\gamma l}\|f\|_2.
\endeq
The reason of reducing to the above estimate is that its left hand side is a linear operator, which makes it possible to use the Fourier method. Moreover, notice that the function $u$ does not depend on the variable $y$, hence by applying Plancherel's theorem in $y$, it suffices to prove 
\beq
\left\|\int_{\R}\hat{f}(x-t, \eta)e^{i u(x)t^2 \eta}\psi_l(u(x)^{1/2}t)\frac{dt}{|t|}\right\|_2 \lesim 2^{-\gamma l}\|f\|_2.
\endeq
Here we use $\hat{f}$ to denote the partial Fourier transform of the function $f$ in the variable $y$. By Fubini's theorem, the last estimate follows simply from
\begin{lem}[Proposition 2.2 in \cite{Guo}]
There exists a universal constant $\gamma>0$ such that for any function $f\in L^2(\R)$ it holds that
\beq
\left\| \int_{\R}f(x-t)e^{iu(x)t^2}\psi_l(u(x)^{1/2}t)\frac{dt}{|t|} \right\|_2 \lesim 2^{-\gamma l}\|f\|_2.
\endeq
\end{lem}
\noindent So far we have finished the proof of the estimate \eqref{0402ee1.1} in Theorem \ref{2604theorem1.1}. $\Box$\\

The estimate \eqref{0504ee1.19} is rather crude, and it only gives the $L^p$ bounds of $\mc{M}^u\circ P_k$ for certain large $p$. Indeed we would like to ask the following 
\begin{question}\label{0504question2.2}
For any $p>1$, does it hold that 
\beq
\|\mc{M}^u\circ P_k (f)\|_p \lesim \|f\|_p,
\endeq
with a constant depending only on $p$? Could we further expect that 
\beq
\|\mc{M}^u f\|_p \lesim \|f\|_p, \forall f\in L^p(\R^2),
\endeq
with a constant again depending only on $p$?
\end{question}

As has been explained in the introduction, the above question has a close connection with the boundedness of the circular maximal function and the Radon variant of the Kakeya problem. Particularly, we hope that by using certain techniques similar to those in \cite{KW}, \cite{Schlag} and \cite{Wolff2}, one can give an affirmative answer to the above question.

\section{Hilbert transform case: Proof of \eqref{0504ee1.12}}
In this section, we will prove the estimate \eqref{0504ee1.12} in Theorem \ref{2604theorem1.1}. That is, for a measurable function $u:\R\to \R$, define 
\beq
\mc{H}^u f(x, y):=\int_{\R}f(x-t, y-u(x)\cdot t^2)\frac{dt}{t},
\endeq
we will show that for all $p>p_0$, where $p_0<2$ is the same as in the last section,  it holds that
\beq\label{2704ee3.2}
\|\mc{H}^u\circ P_k (f)\|_p \lesim \|f\|_p.
\endeq
Here the constant is independent of $k\in \Z$.\\


%
By the same scaling argument as in \eqref{0504ee2.3}, it suffices to consider the case $k=0$. For simplicity, in the rest of the proof, we will always assume that $P_0 f=f.$ The starting point of the proof of the estimate \eqref{2704ee3.2} is the same as the one for the maximal operator case, namely, we first decompose the kernel $1/t$ into a ``high'' frequency part and a ``low'' frequency part. 

To be precise, let $\psi_0:\R\to \R$ be a non-negative smooth function supported on $[-5/2, -1/2]\cup [1/2, 5/2]$ such that 
\beq\label{1602ff2.3}
\sum_{l\in \Z}\psi_l(t)=1, \forall t\neq 0,
\endeq 
where $\psi_l(t):=\psi_0(2^{-l} t)$. This implies 
\beq \label{1602ff2.4}
\sum_{l\in \Z}\psi_l( u(x)^{1/2}t)=1.
\endeq
Hence
\beq
\mc{H}^u f(x, y)=\sum_{l\in \Z}\int_{\R}f(x-t, y- u(x) t^2)\psi_l( u(x)^{1/2}t) \frac{dt}{t}.
\endeq

%
%
%
%
%
%


%
%
%

\noindent {\bf The high frequency part $l\le 0$.} Concerning this part, we show that 
\beq
\sum_{l\le 0}\int_{\R}f(x-t, y- u(x) t^2)\psi_l( u(x)^{1/2}t) \frac{dt}{t}
\endeq
can be bounded by the strong maximal function and the maximal Hilbert transform.\\

Denote by
\beq
\Psi_0(t):=\sum_{l\le 0}\psi_l(t),
\endeq
we obtain
\beq
\begin{split}
& \sum_{l\le 0}\int_{\R} f(x-t, y- u(x) t^2)\psi_l( u(x)^{1/2}t) \frac{dt}{t}\\
&=\int_{\R}f(x-t, y- u(x) t^2)\Psi_0( u(x)^{1/2}t) \frac{dt}{t}.
\end{split}
\endeq
We do the change of variable $u(x)^{1/2}t\to t$
to obtain
\beq
\int_{\R}f(x-u(x)^{-1/2}t, y- t^2)\Psi_0(t) \frac{dt}{t}.
\endeq
The idea is to compare the last expression with 
\beq
\int_{\R}f(x-u(x)^{-1/2}t, y)\Psi_0(t) \frac{dt}{t},
\endeq
hence we write
\beq\label{1602ee2.12}
\begin{split}
& \int_{\R}f(x-u(x)^{-1/2}t, y- t^2)\Psi_0(t) \frac{dt}{t}\\
&=\int_{\R}f(x-u(x)^{-1/2}t, y)\Psi_0(t) \frac{dt}{t}\\
&+\int_{\R}\left(f(x-u(x)^{-1/2}t, y- t^2)-f(x-u(x)^{-1/2}t, y)\right)\Psi_0(t) \frac{dt}{t}.
\end{split}
\endeq
For the former part on the right hand side, we bound it by the one dimensional maximal Hilbert transform. For the latter part, we will be using the fact that the function $f$ has frequency near one in the vertical direction. This will be similar to the estimate \eqref{0504ee2.18}, hence we leave it out.\\
%

\noindent {\bf The low frequency part $l\ge 0$.} By the triangle inequality, it suffices to prove that there exists a constant $\gamma_p>0$ depending only on $p$ such that 
\beq\label{0402ee2.19}
\|\int_{\R}f(x-t, y- u(x) t^2)\psi_l( u(x)^{1/2}t) \frac{dt}{t}\|_p \lesim 2^{-\gamma_p l}\|f\|_p.
\endeq
The case $p=2$ of \eqref{0402ee2.19} has already been proved, see \eqref{1602ff1.20}. To prove \eqref{0402ee2.19} for $p>p_0$, by interpolation, it suffices to prove that 
\beq\label{0402ee2.20}
\|\int_{\R} f(x-t, y- u(x) t^2)\psi_l( u(x)^{1/2}t) \frac{dt}{t}\|_p \lesim \|f\|_p.
\endeq
Notice that as we are not aiming at any exponential decay, the cancellation from the kernel $1/t$ will no longer play any role, hence we simply bound the left hand side of \eqref{0402ee2.20} by 
\beq
\|\sup_{k\in \Z}\mc{M}^{u(x)}_k f(x, y)\|_p,
\endeq
which has been proved to satisfy an $L^p$ bound for all $p> p_0$ in the last section. So far we have finished the proof of the boundedness of the low frequency part, hence the $L^p$ bounds of the operator $\mc{H}^u\circ P_0$ for all $p>p_0$.

\begin{rem}
If we were able to give an affirmative answer to Question \ref{0504question2.2}, then we would also be able to prove the estimate \eqref{0504ee1.12} for all $p>1$.
\end{rem}

Shaoming Guo, Institute of Mathematics, University of Bonn\\
\indent Address: Endenicher Allee 60, 53115, Bonn\\
\indent Email: shaoming@math.uni-bonn.de

\end{document}